# Moment estimates, exponential integrability, concentration inequalities and exit times estimates on evolving manifolds

Robert Baumgarth[1]

[1]Fakultät für Mathematik und Informatik, Augustusplatz 10, 04109 Leipzig, Germany
`robert.baumgarth@math.uni-leipzig.de`

6th February 2024

### Abstract

On a smooth (not necessarily compact) manifold $M$ equipped with a $\mathrm{C}^1$-family of complete Riemannian metrics $g(t)$ and a $\mathrm{C}^{1,\infty}$-family of vector fields $Z(t)$ both indexed by the real interval $[0, T)$ where $T \in (0, \infty]$, we prove moment estimates, exponential integrability, concentration inequalities and exit times estimates for diffusions on complete evolving Riemannian manifolds.

## 1  Introduction

Let $M$ be a smooth (not necessarily compact) $m$-manifold equipped with a $\mathrm{C}^1$-family of complete Riemannian metrics $g(t)$ and a $\mathrm{C}^{1,\infty}$-family of vector fields $Z(t)$ both indexed by the real interval $[0, T)$ where $T \in (0, \infty]$. Let $\nabla^{g(t)}$ and $\Delta_{g(t)}$ be the Levi-Civita connection and the Laplace-Beltrami operator associated to the metric $g(t)$, respectively. If $g(t)$ is independent of $t$, we simply write $\Delta$.

The family of Riemannian metric $g(t)$ describes the evolution of the manifold $M$ under a geometric flow where the time $T$ marks the first time where the curvature possibly blows up. Due to Perelman's breakthrough work making use of the Ricci flow, evolving Riemannian manifolds have been studies intensively over the past decades. Introduced by Hamilton [Ham82] in 1982, Ricci flow describes the evolution of the metric tensor, given some initial metric $g_0$ on $M$, under the nonlinear partial differential equation

$$\partial_t g(t) = -2 \operatorname{Ric}_{g(t)}, \qquad g(0) = g_0,$$

where $\operatorname{Ric}_{g(t)}$ denotes the Ricci curvature of $M$ with respect to the metric $g(t)$. In local harmonic coordinates $x_i$, $\Delta x_i = 0$, then

$$\operatorname{Ric}_{ij} = -\frac{1}{2}\Delta g_{ij} + \text{ lower order terms}$$

acts like a heat equation for metric tensor, reflecting the idea to «smoothen out» the manifold. In this setting Brownian motion with respect to a collection of time-dependent metrics, so called $g(t)$-Brownian motion, was originally defined in 2008 by Arnaudon, Coulibaly-Pasquier and Thalmaier in [ACT08] and later elaborated by Coulibaly-Pasquier in [Cou11]. The theory seamlessly carries over to martingales with respect to a time-dependent connection [GPT15]. Kuwada and Philipowski [KP11] showed that Brownian motion cannot explode in finite time if the metric evolves under backwards super Ricci flow. The key ingredient the proof is the extension of Kendall's [Ken87] well-known Itô formula for the radial process, as the distance function (taken from a fixed reference point) is smooth only outside the cut locus of $M$.

Using the so called Onsager-Machlup functional of an inhomogeneous uniformly elliptic diffusion process on $(M, g(t))$, in [Cou14] Coulibaly-Pasquier showed that the probability that a Brownian motion deviates from a smooth curve by at most a distance $\epsilon > 0$ decays exponentially in $\epsilon$. In Versendaal [Ver21], it is shown how the large deviations are obtained from the large deviations of the (time-dependent) horizontal lift of a $g(t)$-Brownian motion to the frame bundle $\mathcal{O}(M)$ over $M$.

In the recent article [TT20], Thalmaier and Thompson study exponential integrability and exit times of diffusions on sub-Riemannian and metric measure spaces based on [Tho16], where the second author considered the case of $L = \frac{1}{2}\Delta + Z$, where $Z$ is a smooth vector field and $\Delta$ the Laplace-Beltrami operator on a complete Riemannian manifold. In each case the necessary ingredients are an Itô formula and a comparison theorem for the Laplacian.

So let now $M$ be an $m$-manifold and $g(t)$ be a time-dependent family of Riemannian metrics over $M$. Let $\mathbf{L}_t$ be an inhomogeneous uniformly elliptic operator on $M$

$$\mathbf{L}_t = \frac{1}{2}\Delta_{g(t)} + Z(t),$$

where $\Delta_{g(t)}$ is a Laplace Beltrami operator for a metric $g(t)$ and $Z(t, \cdot)$ is a time-dependent vector field over $M$.

Our method is inspired by Thalmaier and Thompson [TT20]: Suppose the diffusion operator $\mathbf{L}_t$ satisfies an inequality of the form $\mathbf{L}_t \varphi \leqslant \nu + \lambda\varphi$ for some constants $\nu, \lambda$ and a suitable function $\varphi$. Then, given an appropriate Itô formula for the corresponding diffusion $X$ on an evolving Riemannian manifold $(M, g(t))$, we can calculate various estimates on the moments of the random variable $\varphi(X_t)$. From those, we derive moment estimates, exponential integrability, concentration inequalities and exit times estimates for the $g(t)$-Brownian motion in the setting of evolving Riemannian geometries.

More precisely, we assume a Lyapunov-like condition: There are constants $\nu \geqslant 1$ and $\lambda \in \mathbb{R}$, such that it holds, for all $0 \leqslant s \leqslant T$,

$$\left(\frac{1}{2}\Delta_{g(s)} + \partial_s\right) r_s^2 \leqslant \nu + \lambda r_s^2,$$

outside the cut locus of $M$, where $\Delta_{g(s)}$ is the evolving Laplace-Beltrami operator with respect to $g(s)$. Here $r_s(x) := \mathrm{d}_s(x, o)$ is the Riemannian distance with respect to $g(s)$ between $x \in M$ and a fixed reference point $o$.

Using an Itô formula for the radial by Kuwada and Philipowski [KP11], we first prove moment estimates for the radial process. By means of the Laguerre polynomials we deduce an exponential estimate, Theorem 3.5,

$$\mathbb{E}^x \mathrm{e}^{\frac{\vartheta}{2} r_t^2(X_t)} \leqslant (1 - \vartheta\Lambda(t))^{-\nu/2} \exp\left(\frac{\vartheta r_t^2(x) \exp^{\lambda t}}{2(1 - \vartheta\Lambda(t))}\right) \tag{1.1}$$

for all $t, \vartheta \geqslant 0$ such that $\vartheta \Lambda(t) < 1$, where $\Lambda(t) := (e^{\lambda t} - 1)/\lambda$. For the geodesic ball $\mathsf{B}^{g(t)}(x_0, r)$ around $x_0 \in M$ with radius $r > 0$, we give the concentration inequality, cf. Theorem 3.6,

$$\lim_{r \to \infty} \frac{1}{r^2} \log \mathbb{P}^x \left( X_t \notin \mathsf{B}^{g(t)}(x_0, r) \right) \leqslant -\frac{1}{2\Lambda(t)}$$

for all $t > 0$ and an exit time estimate in the case of evolving geometries, cf. Theorem 3.8,

$$\mathbb{P}^x \left( \sup_{0 \leqslant s \leqslant t} r_s(X_s) \geqslant r \right) \leqslant (1 - \delta)^{-\frac{\nu}{2}} \exp \left( \frac{r_t^2(x) \delta e^{\lambda t}}{2(1 - \delta)\Lambda(t)} - \frac{\delta r^2}{2\Lambda(t)} \right)$$

for all $t > 0$ and some $\delta \in (0, 1)$.

The outline of the paper is as follows: In § 2, we introduce key notions of Brownian motions with respect to evolving Riemannian manifolds, i.e. so called $g(t)$-Brownian motions. Of particular interest in what follows is an Itô formula for the radial process $r_t(X_t)$ of a $g(t)$-Brownian motion by Kuwada and Philipowski [KP11]. In § 3.1, we derive estimates for the second radial moments using this Itô formula. In §§ 3.2 and 3.3, we will see that higher even radial moments, and an exponential estimate, can be achieved using properties of the Laguerre polynomials for Gaussian random variables, and an exponential identity for a sum of Laguerre polynomials respectively. The concentration inequality follows through an application of the Markov inequality, cf. § 3.4. The exit estimate requires a slight modification of the arguments made before and is proved in § 3.5. In § 4 we show that the argument extends to non-symmetric diffusion by a slight modification of the assumption (A) and give conditions when the modified condition (B) is satisfied.

Let us close the introduction with some conventions on general notation used. The metric $g(t)$ on $M$ is denoted by $\langle \cdot, \cdot \rangle_{g(t)} := g(t)(\cdot, \cdot)$ interchangeably. For any two-tensor $T_t$ and any function $f$ (or constant $K$), we write $T_t \geqslant f$ (or $T_t \geqslant K$ respectively) if $T_t(X, X) \geqslant f \langle X, X \rangle_t = f \left| X \right|_t$ holds for $X \in \Gamma(\mathsf{T}M)$. Here $\Gamma(\mathsf{T}M)$ denotes the smooth sections of the tangent bundle over $M$. By $\mathbb{E}^{(s,x)}$ and $\mathbb{P}^{(s,x)}$ denote, respectively, the expectation and the probability of the underlying process starting from $x$ at time $s$. When $s = 0$, we simply write $\mathbb{E}^x := \mathbb{E}^{(0,x)}$ and $\mathbb{P}^x = \mathbb{P}^{(0,x)}$. Throughout the article, we will use the Einstein summation convention: If an index appears twice, once as a subscript and once as a superscript, we skip the summation symbol, e.g. we write $g(ue_\alpha, ue_\beta)V^{\alpha\beta}$ instead of $\sum_{\alpha, \beta \leqslant m} g(ue_\alpha, ue_\beta)V^{\alpha\beta}$, where $m := \dim M$ is the dimension of $M$.

## 2  Brownian Motion on Evolving Manifolds and Itô Formula

Let $(M, g(t))_{t \in I}$ be an evolving manifold indexed by $I = [0, T)$ for some $T \in (0, \infty]$. Let $\nabla^{g(t)}$ the Levi-Civita connection with respect to $g(t)$. We consider the bundle

$$\mathsf{T}M \overset{\pi}{\longrightarrow} \mathbb{M},$$

where $\pi$ is the projection onto space-time $\mathbb{M} := M \times I$. It is due to Hamilton [Ham93] that there is a natural space-time connection $\nabla$ on $\mathsf{T}M$ considered as a bundle over space-time $\mathbb{M}$ given by

$$\begin{aligned}
\nabla_v X &= \nabla_v^{g(t)} X, \\
\nabla_{\partial_t} X &= \partial_t X + \frac{1}{2} \partial_t g_t(X, \cdot)^{\sharp g_t},
\end{aligned} \tag{2.1}$$

for all $v \in (T_x M, g(t))$ and all time-dependent vector fields $X$ on $M$. In particular, this connection is compatible with the metric, i.e.

$$\frac{\mathrm{d}}{\mathrm{d}t} |X|_{g(t)}^2 = 2 \left\langle X, \nabla^{g(t)} X \right\rangle_{g(t)}.$$

The construction presented in the remainder of this section closely follows the classical ideas for the static case, see e.g. [Tha16; HT94] for the static case and [ACT08; Cou11; CT23; Che17] for the extension to the evolving setting.

**Remark 2.1.** Let $G = O(m)$, with $\dim M = m$, and consider the G-principal bundle $P := \mathcal{O}(M) \xrightarrow{\pi} \mathbb{M}$ of orthonormal frames with fibres

$$\mathcal{O}_{(x,t)}(M) = \left\{ u : \mathbb{R}^m \to (T_x M, g(t)) : u \text{ isometry} \right\},$$

where $a \in G$ acts on $\mathcal{O}(M)$ from the right

$$u \triangleleft g : \mathbb{R}^m \xrightarrow{g} \mathbb{R}^m \xrightarrow{u} T_x M,$$

The connection $\nabla^{g(t)}$ gives rise to a G-invariant splitting of the sequence

$$0 \longrightarrow \ker \mathbf{d}\pi \longrightarrow T\mathcal{O}(M) \underset{h}{\overset{\mathbf{d}\pi}{\longrightarrow}} \pi^* T\mathbb{M} \longrightarrow 0$$

which induces a decomposition of $T\mathcal{O}(M)$,

$$T\mathcal{O}(M) = V \oplus H := \ker \mathbf{d}\pi \oplus h(\pi^* T\mathbb{M}).$$

For each $u \in \mathcal{O}(M)$, the **horizontal space** $H_u$ **at** $u$ is constituted via G-invariance, i.e. it holds $H_{u \triangleleft g} \mathcal{O}(M) = (\triangleleft g)_* H_u \mathcal{O}(M)$, for the G-right action $\triangleleft g$ on $M$ and the **vertical space** $V_u$ **at** $u$ is given by $V_u = \left\{ v \in T_u \mathcal{O}(M) : (\mathbf{d}\pi)v = 0 \right\}$. The bundle isomorphism

$$h : \pi^* T\mathbb{M} \xrightarrow{\sim} H \to T\mathcal{O}(M), \qquad h_u : T_{\pi(u)} \mathbb{M} \xrightarrow{\sim} H_u, \qquad u \in \mathcal{O}(M),$$

is called the **horizontal lift** of the G-connection (cf. fig. 0.1 below).

**Corollary 2.2.** *To each $X + b\partial_t \in T_{(x,t)} \mathbb{M}$ and each frame $u \in \mathcal{O}_{(x,t)}(M)$, there is a unique horizontal lift $X^\uparrow + b\partial_t^\uparrow \in H_u$ of $X + b\partial_t$ such that*

$$\pi_* \left( X^\uparrow + b\partial_t^\uparrow \right) = X + b\partial_t,$$

*where $X^\uparrow$ is the horizontal lift of $X$ with respect to the metric $g(t)$, and $\partial_t^\uparrow = \frac{\mathrm{d}}{\mathrm{d}s}\big|_{s=0} u_s$ with $u_s$ is the horizontal lift based at $u$ of the curve $s \mapsto (x, t+s)$.*

**Remark 2.3.** Qua horizontal lift, there are $m$ well-defined unique (standard) **horizontal vector fields**

$$H_i^{g(t)} \in \Gamma(T\mathcal{O}(M)), \qquad H_i^{g(t)}(u) := (ue_i)^\uparrow = h_u^t(ue_i), \qquad i = 1, \dots, m,$$

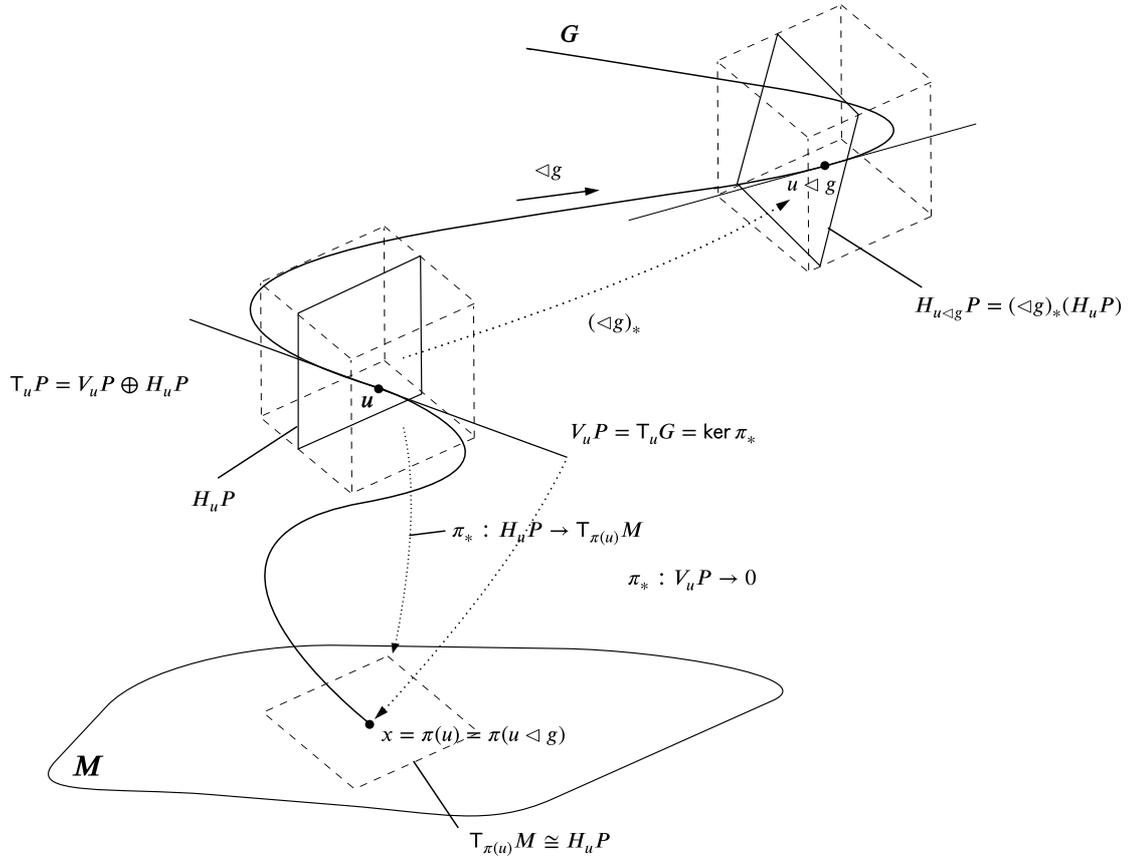

Figure 0.1: *G*-connection of a principal *G*-bundle

i.e. whose projection is the $i^{\text{th}}$ unit vector $ue_i$ of the orthonormal frame, $\pi_* H_i^{g(t)}(u) = ue_i$, where $(e_i)_{1 \leqslant i \leqslant m}$ is the canonical basis for $\mathbb{R}^m$. Those vector fields give rise to a Laplacian of Hörmander type on $\mathcal{O}(M)$, i.e. that can be written as sums of squares as $\mathcal{O}(M)$ is parallelizable: We call

$$\Delta_{\mathcal{O}(M)} := \sum_{i=1}^{n} (H_i^{g(t)})^2$$

**Bochner's horizontal Laplacian on $\mathcal{O}(M)$.** Note that by the fundamental relation $\Delta_{\mathcal{O}(M)} \pi^* = \pi^* \Delta$ (cf. [Mal97, p. 286]), for every $f \in C^\infty(M)$ and $\tilde{f} := f \circ \pi$ its lift to $\mathcal{O}(M)$, hence

$$\Delta_{\mathcal{O}(M)} \tilde{f}(u) = \Delta f(x) \quad \text{for any } u \in \mathcal{O}(M) \text{ with } x = \pi(u).$$

Finally, the horizontal lift of a vector field $Z(t) \in \Gamma(\mathsf{T}M)$ to $\mathscr{F}$ is given by

$$H_{Z(t)}^t \in \Gamma(\mathsf{T}\mathscr{F}), \qquad H_{Z(t)}^t(u) = h_u(\nabla^{g(t)} Z(t, \cdot)_x), \qquad \pi(u) = (x, t).$$

Before stating the Stratonovich SDE on the frame bundle $\mathscr{F}$, we give the explicit formula for horizontal lift $\partial_t^\dagger$ of the time direction $\partial_t$ in Hamilton's connection (2.1), cf. e.g. [Cou11, Proposition 1.2.].

**Proposition 2.4.** *For any $u \in \mathscr{F}$ the horizontal lift $\partial_t^\uparrow$ of $\partial_t$ with respect to the Hamilton's connection $\nabla$ (cf. (2.1)) is given by*

$$\partial_t^\uparrow(u) = -\frac{1}{2}\partial_t g_t(ue_\alpha, ue_\beta)V^{\alpha\beta}(u),$$

*where $(e_1, \ldots, e_m)$ denotes the canonical basis for $\mathbb{R}^m$ and $(V^{\alpha\beta}(u))_{\alpha,\beta=1,\ldots,m}$ is the canonical basis of vertical vector fields over $\mathscr{F}$.*

Next, we consider the Stratonovich SDE on $\mathscr{F}$:

$$\mathrm{d}U = H_i^{g(t)}(U) \circ \mathrm{d}B^i + H_{Z_t}^{g(t)}(U)\mathrm{d}t - \frac{1}{2}\partial_t g(t)(Ue_\alpha, Ue_\beta)V^{\alpha\beta}(U)\mathrm{d}t$$
$$U_s \in \mathscr{F}_s(M), \qquad t \in [0, T), \tag{2.2}$$

where $B^i$ denotes a standard Brownian motion on $\mathbb{R}^m$ with generator $\frac{1}{2}\Delta_{\mathbb{R}^m}$. Since the horizontal vector fields $H_{Z(t)}^{g(t)}$ are in $C^{1,\infty}$, the SDE (2.2) has a unique solution up to its lifetime $\zeta := \lim_n \zeta_n$ with

$$\zeta_n := \inf\left\{t \in [s, T) : \mathrm{d}_s(\pi(U_s), \pi(U_t)) \geqslant n\right\}, \qquad n \in \mathbb{N}, \qquad \inf \varnothing := T,$$

and $\mathrm{d}_t$ the Riemannian distance with respect to $g(t)$.

Note that (2.2) has local Lipschitz coefficients and hence a unique strong solution with continuous paths. By Proposition 2.4, the choice is also canonical and forces the solution to the SDE to be a $g(t)$-isometry, cf. [Cou11].

If $U$ is a solution to the SDE (2.2) then

$$\pi(U_t) = (X_t, t),$$

where $X$ is a diffusion process on $M$ generated by $\mathsf{L} = \frac{1}{2}\Delta_{g(t)} + Z(t)$, is called $g(t)$-**Brownian motion with drift** $Z(t)$ **on** $\mathbb{M}$. If $Z = 0$, then $X$ is called a $g(t)$-**Brownian motion on** $\mathbb{M}$.

The induced **parallel transport along** $X_t$ is now defined as

$$/\!/_{r,t} := U_t U_r^{-1} : (\mathrm{T}_{X_r}M, g(r)) \to (\mathrm{T}_{X_t}M, g(t)), \qquad s \leqslant r \leqslant t < T.$$

If $r = 0$, we simply write $/\!/_t := /\!/_{0,t}$. Note that, if $g(t)$ is independent of $t$, then (2.2) reduces to the usual parallel transport and that the $/\!/_{r,t}$ are, by definition, isometries.

More precisely [CT23, Proposition 2.5], it holds, for any $C^2$-function $F : \mathscr{F} \to \mathbb{R}$,

$$\mathrm{d}F(U) = H_i^{g(t)}F(U)\mathrm{d}B^i + \Delta_{\mathscr{O}(M)}F(U)\mathrm{d}t + H_{Z_t}^{g(t)}(U)\mathrm{d}t + \partial_t^\uparrow F(U)\mathrm{d}t$$

and it holds, for any $C^2$-function $f : \mathbb{M} \to \mathbb{R}$,

$$\mathrm{d}f(X) = \partial_t f(X)\mathrm{d}t + (Ue_i f)\mathrm{d}B^i + \mathsf{L}_t f(X)\mathrm{d}t.$$

If $X_s = x$ is common to denote by

$$X_t = X_t^{(x,s)}, \qquad t \geqslant s,$$

which solves the equation

$$\mathrm{d}X_t^{(s,x)} = U_t \circ \mathrm{d}B_t + Z_t(X_t^{(s,x)})\mathrm{d}t, \qquad X_s^{(s,x)} = x.$$

By Itô's formula, for any $f \in C^2(M)$,

$$f(X_t^{(s,x)}) - f(x) - \int_s^t \mathsf{L}_r f(X_r^{(s,x)}) \mathrm{d}r = \int_s^t \left\langle \nabla^r f(X_r^{(s,x)}), /\!/_s \mathrm{d}B_r \right\rangle_s, \qquad s \leqslant t < T,$$

is a martingale up to the lifetime $\zeta$. In the case $s = 0$, we abbreviate $X_t^{(0,x)}$ by $X_t^x$.

Let $r_t(x) := r(t, x) := \mathrm{d}_t(x, o)$ denote the Riemannian distance with respect to $g(t)$ between $x$ and a fixed reference point $o$. We recall that the cut locus $\mathrm{Cut}_{g(t)}(o)$ of $o$ for the distance $\mathrm{d}_t$ is defined as the complement of the set of points $y$ in M for which there exists a unique length minimising geodesic connecting $o$ with $y$ and such that $o$ and $y$ are not conjugate. The global cut locus $\mathrm{Cut}_{g(t)}(M)$ is defined by

$$\mathrm{Cut}_{g(t)}(M) := \left\{ (x, y) \times M \times M \,:\, y \in \mathrm{Cut}_{g(t)}(x) \right\}.$$

It is well-known that the set $M \smallsetminus \mathrm{Cut}_{g(t)}(o)$ is a closed subset in $M$ and that the function $\mathrm{d}_t^2$ is smooth on $(M \times M) \smallsetminus \mathrm{Cut}_{g(t)}(M)$.

A key tool in the proof of our main result will be the following Itô formula for the radial process $r_t(X_t)$ taken from [KP11], extending the well-known result of Kendall [Ken87] for a static Riemannian metric. Note that the usual Itô formula cannot be applied because the distance function is not smooth at the cut locus.

**Theorem 2.5** ([KP11, Theorem 2]). *There is a non-decreasing continuous process $\ell$ which increases only when $X_t \in \mathrm{Cut}_{g(t)}(o)$ such that*

$$r_t(X_t) = r_0(X_0) + \int_0^t \langle \nabla r_s(X_s), /\!/_s \mathrm{d}B_s \rangle + \int_0^t \left( \frac{1}{2} \Delta_{g(s)} r_s(X_s) + \partial_s r_s(X_s) \right) \mathrm{d}s - \ell_t. \tag{2.3}$$

In particular, the quadratic variation $\langle \beta \rangle_t$ of the martingale term

$$\beta_t := \int_0^t \langle \nabla r_s(X_s), /\!/_s \mathrm{d}B_s \rangle$$

is given by

$$\begin{aligned}
\langle \beta \rangle_t &= \sum_{i=1}^m \int_0^t \left[ (/\!/_s e_i) r_s(X_s) \right]^2 \mathrm{d}s \\
&= \int_0^t \left| \nabla^{g(s)} r_s(X_s) \right|^2 \mathrm{d}s \\
&= t.
\end{aligned}$$

**Remark 2.6.** Note that, for all $0 \leqslant t \leqslant T$, the Lebesgue measure of the set of times when $X_t \in \mathrm{Cut}_{g(t)}(o)$ has measure zero almost surely [KP11, Lemma 3] and the function $(t, x) \mapsto r_t(x)$ is smooth whenever $x \in \{o\} \cup \mathrm{Cut}_{g(t)}(o)$ [KP11, Lemma 4].

# 3  Symmetric Diffusions $Z = 0$

We start with the case of a symmetric diffusion, i.e. $Z = 0$ under the following

**Assumption 3.1.** There are constants $\nu \geqslant 1$ and $\lambda \in \mathbb{R}$, such that it holds, for all $0 \leqslant s \leqslant T$,

$$\left(\frac{1}{2}\Delta_{g(s)} + \partial_s\right) r_s^2 \leqslant \nu + \lambda r_s^2, \tag{A}$$

on $M \setminus \mathrm{Cut}_{g(s)}(o)$ for the evolving Laplace-Beltrami operator $\Delta_{g(s)}$ with respect to $g(s)$. Here $r_s(x) := \mathrm{d}_s(x, o)$ is the Riemannian distance with respect to $g(s)$ between $x \in M$ and a fixed reference point $o$.

### 3.1 Second Radial Moments

We start by deriving an estimate for the second moment of the radial function $r_t(X_t)$. In the next section, we will show that higher even radial moments can be deduced by an induction argument using Laguerre polynomials.

**Theorem 3.2.** *Assume that* (A) *holds. Then* $X$ *is nonexplosive and*

$$\mathbb{E}^x r_t^2(X_t) \leqslant r_t^2(x)\mathrm{e}^{\lambda t} + \nu \Lambda(t),$$

*where*

$$\Lambda(t) := \frac{\mathrm{e}^{\lambda t} - 1}{\lambda} \qquad \text{for all } t \geqslant 0. \tag{3.1}$$

**Proof.** Let $(D_i)_{i=1}^{\infty}$ be an exhaustion of $M$ by regular domains and denote by $\tau_{D_i}$ the first exit time of $X(x)$ from $D_i$. Note that $\tau_{D_i} < \tau_{D_{i+1}}$ and that this sequence of stopping times announces the explosion time $\zeta(x)$. Then by the relation

$$\frac{1}{2}\Delta_{g(t)}r_t^2 = r_t\Delta_{g(t)}r_t + \left|\nabla^{g(t)}r\right|_{g(t)}^2$$

and the Itô's product rule applied to (2.3), it follows that

$$r_t^2(X_{t\wedge\tau_{D_i}}) = r_t^2(x) + \int_0^{t\wedge\tau_{D_i}} \left(\frac{1}{2}\Delta_{g(s)}r_s(X_s) + \partial_s r_s(X_s)\right)\mathrm{d}s - 2\int_0^{t\wedge\tau_{D_i}} r_s(X_s)(\mathrm{d}\beta_s - \mathrm{d}\ell_s).$$

holds $\mathbb{P}^x$-almost surely. Since the domains $D_i$ are of compact closure $\beta$ is a martingale. Hence,

$$\mathbb{E}^x r_t^2(X_{t\wedge\tau_{D_i}}) = r_t^2(x) + \mathbb{E}^x \int_0^{t\wedge\tau_{D_i}} \left(\frac{1}{2}\Delta_{g(s)}r_s(X_s) + \partial_s r_s(X_s)\right)\mathrm{d}s - 2\mathbb{E}^x \int_0^{t\wedge\tau_{D_i}} r_t(X_s)\mathrm{d}\ell_s.$$

for all $t \geqslant 0$. Next, we want to apply Gronwall's inequality. But note that the coefficient $\lambda$ is allowed to be negative. Therefore, we split

$$\mathbb{E}^x r_t^2(X_{t\wedge\tau_{D_i}}) = \mathbb{E}^x \left(\mathbb{1}_{\{t<\tau_{D_i}\}} r_t^2(X_t)\right) + \mathbb{E}^x \left(\mathbb{1}_{\{t\geqslant\tau_{D_i}\}} r_t^2(X_{\tau_{D_i}})\right)$$

and note that the two functions

$$t \mapsto \mathbb{E}^x \int_0^{t\wedge\tau_{D_i}} r_t(X_s)\mathrm{d}\ell_s \qquad \text{and} \qquad t \mapsto \mathbb{E}^x \left(\mathbb{1}_{\{t\geqslant\tau_{D_i}\}} r_t^2(X_{\tau_{D_i}})\right)$$

are non-decreasing. Define a function $f_x^{i,2}(t)$ by

$$f_x^{i,2}(t) := \mathbb{E}^x \left(\mathbb{1}_{\{t<\tau_{D_i}\}} r_t^2(X_t)\right).$$

Then $f_x^{i,2}$ is differentiable and by assumption we have the differential inequality

$$\frac{\mathrm{d}}{\mathrm{d}t} f_x^{i,2}(t) \leqslant \nu + \lambda f_x^{i,2}(t)$$
$$f_x^{i,2}(0) = r_t^2(x)$$

for all $t \geqslant 0$. By Gronwall's inequality, we get

$$\mathbb{E}^x \left( \mathbb{1}_{\{t < \tau_{D_i}\}} r_t^2(X_t) \right) \leqslant r_t^2(x) \mathrm{e}^{\lambda t} + \nu \frac{\mathrm{e}^{\lambda t} - 1}{\lambda} \tag{3.2}$$

for all $t \geqslant 0$. Choosing $D_i = \mathrm{B}(i, o) := \left\{ y \in M : r_t(y) < i \right\}$, the inequality (3.2) yields

$$\mathbb{P} \left( \tau_{\mathrm{B}(i,o)} \leqslant t \right) \leqslant \frac{r_t^2(x) \mathrm{e}^{\lambda t} + \nu \Lambda(t)}{i^2}$$

for all $t \geqslant 0$, which implies that X is nonexplosive. Hence, the assumption follows by Beppi Levi's theorem. ∎

A simply application of Jensen's inequality yields the following

**Corollary 3.3.** *Assume that* (A) *holds. Then*

$$\mathbb{E}^x r_t(X_t) \leqslant \sqrt{r_t(x) \mathrm{e}^{\lambda t} + \nu \Lambda(t)}.$$

## 3.2 Higher Even Radial Moments

In this section we calculate higher even radial moments by means of the **Laguerre polynomials** $L_p^\alpha(z)$ defined, for $p \in \mathbb{N}_0$ and $\alpha > -1$, by

$$L_p^\alpha(z) := \mathrm{e}^z \frac{z^{-\alpha}}{p!} \frac{\partial^p}{\partial z^p} \left( \mathrm{e}^{-z} z^{p+\alpha} \right). \tag{3.3}$$

Then for any real-valued Gaussian random variable with mean $\mu$ and variance $\sigma^2$ (cf. e.g. [Leb72])

$$\mathbb{E} Y^{2p} = (2\sigma^2)^p p! L_p^{-1/2} \left( -\frac{\mu^2}{2\sigma^2} \right).$$

In particular, if $X$ is a standard Brownian motion on $\mathbb{R}$, then

$$\mathbb{E}^x \left| X_t \right|^{2p} = (2t)^p p! L_p^{-1/2} \left( -\frac{|x|^2}{2t} \right) \qquad \text{for all } t \geqslant 0.$$

**Theorem 3.4.** *Assume that* (A) *holds. Then, for* $p \in \mathbb{N}$,

$$\mathbb{E}^x r_t^{2p}(X_t) \leqslant (2\Lambda(t))^p p! L_p^{\nu/2-1} \left( -\frac{r_t^2(x) \mathrm{e}^{\lambda t}}{2\Lambda(t)} \right)$$

*for all* $t \geqslant 0$, *where* $\Lambda(t)$ *is defined as in Theorem 3.2.*

**Proof.** By (A) it easily follows that, on $M \smallsetminus \mathrm{Cut}_{g(t)}(o)$ and for $p \in \mathbb{N}$,

$$\left( \frac{1}{2} \Delta_{g(s)} + \partial_s \right) r^{2p} \leqslant p(\nu + 2(p-1)) r^{2p-2} + p\lambda r^{2p}.$$

Let $\tau_{D_i}$ defined as in the proof of Theorem 3.2. By Itô's formula, we have

$$r_t^{2p}(X_{t \wedge \tau D_i}) = r_0^{2p}(x) + 2p \int_0^{t \wedge \tau_{D_i}} r_s^{2p-1}(X_s)\left(\mathrm{d}\beta_s - \mathrm{d}l_s\right) + \int_0^{t \wedge \tau_{D_i}} \left(\frac{1}{2}\Delta_{g(s)}r_s(X_s) + \partial_s r_s(X_s)\right)\mathrm{d}s$$

for all $t \geqslant 0$, almost surely. So, similarly, if we define functions

$$f_x^{i,2p}(t) := \mathbb{E}^x\left(\mathbb{1}_{\{t < \tau_{D_i}\}}r_t^{2p}(X_t)\right).$$

and have differential inequalities

$$\frac{\mathrm{d}}{\mathrm{d}t}f_x^{i,2p}(t) \leqslant p(v + 2(p-1))f_x^{i,2(p-1)}(t) + p\lambda f_x^{i,2p}(t) \tag{3.4}$$

$$f_x^{i,2p}(0) = r_t^{2p}(X) \tag{3.5}$$

By Gronwall's inequality, we get

$$f_x^{i,2p}(t) \leqslant \left(r^{2p}(x) + p(v+2(p-1))\int_0^t f_x^{i,2(p-1)}(s)\mathrm{e}^{-p\lambda s}\mathrm{d}s\right)\mathrm{e}^{p\lambda t} \tag{3.6}$$

for all $t \geqslant 0$ and $p \in \mathbb{N}$. The next step in the proof is to use induction on $p$ to show that

$$f_x^{i,2p}(t) \leqslant \sum_{k=0}^p \binom{p}{k}(2\Lambda(t))^{p-k}r^{2k}(x)\frac{\Gamma(\frac{v}{2}+p)}{\Gamma(\frac{v}{2}+k)}\mathrm{e}^{p\lambda t} \tag{3.7}$$

for all $t \geqslant 0$ and $p \in \mathbb{N}$, where $\Lambda$ is defined by (3.1). The case $p = 1$ is covered by (3.2) in the proof of Theorem 3.2. Suppose that the inequality holds for some $p - 1$, then by (3.6), we have

$$f_x^{i,2p}(t) \leqslant \left(r^{2p}(t,x) + p(v+2(p-1))\sum_{k=0}^{p-1}\binom{p-1}{k}r^{2k}(t,x)\frac{\Gamma(\frac{v}{2}+p-1)}{\Gamma(\frac{v}{2}+k)}\tilde{\lambda}(t)\right)\mathrm{e}^{p\lambda t} \tag{3.8}$$

for all $t \geqslant 0$, where $\tilde{\lambda}(t) = \int_0^t (2\lambda(s))^{p-1-k}\mathrm{e}^{-\lambda s}\mathrm{d}s$. Using $2(p-k)\tilde{\lambda}(t) = (2\lambda(t))^{p-k}$ and properties of the Gamma function it is straightforward to deduce inequality (3.7) from inequality (3.8), which completes the inductive argument. Since $v \geqslant 1$ we can then apply the relation

$$L_p^\alpha(z) = \sum_{k=0}^p \frac{\Gamma(p+\alpha+1)}{\Gamma(k+\alpha+1)}\frac{(-z)^k}{k!(p-k)!}, \tag{3.9}$$

which can be proved using Leibniz' formula, to see that

$$\sum_{k=0}^p \binom{p}{k}(2\lambda(t))^{p-k}r^{2k}(x)\frac{\Gamma(\frac{v}{2}+p)}{\Gamma(\frac{v}{2}+k)} = (2\Lambda(t)\mathrm{e}^{\lambda t})^p p!L_p^{\frac{v}{2}-1}\left(-\frac{r_t^2(x)\mathrm{e}^{\lambda t}}{2\Lambda(t)}\right)$$

and so by inequality (3.7) it follows that

$$f_x^{i,2p}(t) \leqslant (2\Lambda(t))^p p!L_p^{\frac{v}{2}-1}\left(-\frac{r_t^2(x)\mathrm{e}^{\lambda t}}{2\Lambda(t)}\right) \tag{3.10}$$

for $t \geqslant 0$ and $i, p \in \mathbb{N}$. The result follows from this by Beppo Levi's theorem. $\blacksquare$

### 3.3 Exponential estimate

For $|\gamma| < 1$ the Laguerre polynomials also satisfy the identity

$$\sum_{p=0}^{\infty} \gamma^p L_p^\alpha(z) = (1-\gamma)^{-(\alpha+1)} e^{-\frac{z\gamma}{1-\gamma}}. \tag{3.11}$$

Using this identity and (3.3) from the previous section, for any real-valued Gaussian random variable $Y$ with mean $\mu$ and variance $\sigma^2$ it follows that, for $\vartheta \geqslant 0$ with $\vartheta\sigma^2 < 1$,

$$\mathbb{E} e^{\frac{\vartheta}{2}|Y|^2} = \left(1 - \vartheta\sigma^2\right)^{-1/2} \exp\left(\frac{\vartheta|\mu|^2}{2(1-\vartheta\sigma^2)}\right).$$

In particular, if $X$ is a standard Brownian motion on $\mathbb{R}$ starting from $x \in \mathbb{R}$, then for $t \geqslant 0$ we have that

$$\mathbb{E} e^{\frac{\vartheta}{2}|X_t(x)|^2} = (1-\vartheta t)^{-1/2} \exp\left(\frac{\vartheta|x|^2}{2(1-\vartheta t)}\right).$$

so long as $\vartheta t < 1$.

We use those identities to prove the following exponential estimate in the case of Brownian motion on a evolving Riemannian manifold.

**Theorem 3.5.** *Assume that* (A) *holds. Then*

$$\mathbb{E}^x e^{\frac{\vartheta}{2} r_t^2(X_t)} \leqslant (1 - \vartheta\Lambda(t))^{-\nu/2} \exp\left(\frac{\vartheta r_t^2(x) e^{\lambda t}}{2(1-\vartheta\Lambda(t))}\right) \tag{3.12}$$

*for all $t, \vartheta \geqslant 0$ such that $\vartheta\Lambda(t) < 1$, where $\Lambda(t)$ is defined as in Theorem 3.2.*

**Proof.** By inequality (3.10) and the identity for the Laguerre polynomials (3.11)

$$\begin{aligned}
\mathbb{E}^x\left(\mathbb{1}_{\{t < \tau_{D_i}\}} e^{\frac{\vartheta}{2} r_t^2(X_t)}\right) &= \sum_{p=0}^{\infty} \frac{\vartheta^p}{2^p p!} f_x^{i,2p}(t) \\
&\leqslant \sum_{p=0}^{\infty} (\vartheta t \lambda(t))^p L_p^{\frac{\nu}{2}-1}\left(-\frac{r_t^2(x) e^{\lambda t}}{2\Lambda(t)}\right) \\
&= (1 - \vartheta\Lambda(t))^{-\frac{\nu}{2}} \exp\left(\frac{\vartheta r_t^2(x) e^{\lambda t}}{2(1-\vartheta\Lambda(t))}\right)
\end{aligned}$$

where we justify switching the order of integration with the stopping time. The result follows by Beppo Levi's theorem. ∎

### 3.4 Concentration inequality

It is well-known that for $X$ a Brownian motion on $\mathbb{R}$ starting at $x \in \mathbb{R}$

$$\lim_{r \to \infty} \frac{1}{r^2} \log \mathbb{P}^x\left(X_t \notin B^{\mathbb{R}}(x_0, r)\right) = -\frac{1}{2t} \tag{3.13}$$

for all $t > 0$ and for the Euclidean open ball around $x_0 \in \mathbb{R}$ with radius $r > 0$.

In particular, note that the righthand side is independent of the dimension $m$. For a Brownian motion $X$ on $M$ the relation can be expected to still hold true as a in inequality by a heat kernel comparison argument under suitable curvature conditions. In fact, an asymptotic estimate of this form it is shown in [Str00, Theorem 8.62] if the Ricci curvature is bounded below and was later proved in [Tho16] studying the distance between a Brownian motion and a closed embedded submanifold of a complete Riemannian manifold.

Now, let $\mathsf{B}^{g(t)}(x_0, r)$ be the open geodesic ball with respect to the metric $g(t)$ around $x_0 \in M$ with radius $r > 0$.

**Theorem 3.6.** *Assume that* (A) *holds. Then*

$$\lim_{r \to \infty} \frac{1}{r^2} \log \mathbb{P}^x \left( X_t \notin \mathsf{B}^{g(t)}(x_0, r) \right) \leqslant -\frac{1}{2\Lambda(t)}$$

*for all $t > 0$, where $\Lambda(t)$ is defined as in Theorem 3.2.*

**Proof.** For $\vartheta \geqslant 0$ and $r > 0$, by the Markov inequality and Theorem 3.5

$$
\begin{aligned}
\mathbb{P}^x \left( X_t \notin \mathsf{B}^{g(t)}(x_0, r) \right) &= \mathbb{P}^x \left( r_t(X_t) \geqslant r \right) \\
&= \mathbb{P}^x \left( \mathrm{e}^{\frac{\vartheta}{2} r_t^2(X_t)} \geqslant \mathrm{e}^{\frac{\vartheta}{2} r^2} \right) \\
&\leqslant \mathrm{e}^{-\frac{\vartheta}{2} r^2} \mathbb{E}^x \mathrm{e}^{\frac{\vartheta}{2} r_t^2(X_t)} \\
&\leqslant \left( 1 - \vartheta \Lambda(t) \mathrm{e}^{\lambda t} \right)^{-\frac{\nu}{2}} \exp \left( \frac{\vartheta r_t^2(x) \mathrm{e}^{\lambda t}}{2(1 - \vartheta \Lambda(t) \mathrm{e}^{\lambda t})} - \frac{\vartheta r^2}{2} \right)
\end{aligned}
$$

so long as $\vartheta \Lambda(t) \mathrm{e}^{\lambda t} < 1$. If $t > 0$ then choosing $\vartheta = \delta(\Lambda(t) \mathrm{e}^{\lambda t})^{-1}$ shows that for any $\delta \in [0, 1)$ and $r > 0$ we have the estimate

$$\mathbb{P}^x \left( X_t \notin \mathsf{B}(x_0, r) \right) \leqslant (1 - \delta)^{-\frac{\nu}{2}} \exp \left( \frac{r_t^2(x) \delta \mathrm{e}^{\lambda t}}{2(1 - \delta) t \Lambda(t)} - \frac{\delta r^2}{2\Lambda(t)} \right) \tag{3.14}$$

But $\delta$ can be chosen arbitrarily close to 1, so the assertion follows after taking the limit. ∎

**Remark 3.7.** Note that, by the choice of $\Lambda(t)$, we have $\Lambda(t) \xrightarrow{\lambda \to 0} 1$, so Theorem 3.6 reflects the well-known result (3.13) in the Euclidean case.

## 3.5 Exit time estimate

For the case $\lambda \geqslant 0$, we get the exit time estimate:

**Theorem 3.8.** *Fix $r > 0$ and assume that* (A) *holds on the ball* $\mathsf{B}^{g(t)}(x_0, r)$. *Then*

$$\mathbb{P}^x \left( \sup_{0 \leqslant s \leqslant t} r_s(X_s) \geqslant r \right) \leqslant (1 - \delta)^{-\frac{\nu}{2}} \exp \left( \frac{r_t^2(x) \delta \mathrm{e}^{\lambda t}}{2(1 - \delta)\Lambda(t)} - \frac{\delta r^2}{2\Lambda(t)} \right)$$

*for all $t > 0$ and $\delta \in (0, 1)$.*

**Proof.** The proof requires a slight modification of the argument we used to derive Theorem 3.5 and Theorem 3.6.

We denote by $\tau_r$ the first exit time of $X$ from the open geodesic ball. Applying Itô formula as in the proof of Theorem 3.4, using the assumption $\lambda \geqslant 0$, we obtain the slightly different estimate

$$\mathbb{E}^x r_t^{2p}(X_{t \wedge \tau_r}) \leqslant r_t^{2p}(x) + \frac{p}{2}\left(\nu + 2(p-1)\right) \int_0^t \mathbb{E}^x r_t^{2p-2}(X_{s \wedge \tau_r}) \mathrm{d}s$$
$$+ \frac{p\lambda}{2} \int_0^t \mathbb{E}^x r_t^{2p}(X_{s \wedge \tau_r}) \mathrm{d}s.$$

Following the lines of the inductive argument as before, we deduce moment estimates which can then be summed, as in 3.2, to obtain an exponential estimate for the stopped process. Finally we estimate

$$\mathrm{e}^{\frac{\vartheta}{2}r^2}\mathbb{P}^x\left(\sup_{0 \leqslant s \leqslant t} r_s(X_s) \geqslant r\right) \leqslant \mathbb{E}^x \mathrm{e}^{\frac{\vartheta}{2}r_t^2(X_{t \wedge \tau_r})}$$

and the the right-hand side of the inequality is bounded by the right-hand side of (3.12). Choosing $\vartheta$ as in the proof of Theorem 3.6 yields the claim. ∎

## 4 Generalisation to Non-symmetric Diffusions

All results carry over almost verbatim to the more general setting of non-symmetric diffusions. In [KP11], non-explosion is shown under an additional assumption on the covariant derivative of $Z(t)$ and an additional process given by the Feller test for explosion of SDEs. In our case, we slightly modify our assumption (A).

So from now on, let $X_t$ be a time-dependent diffusion with generator

$$\mathsf{L}_t = \frac{1}{2}\Delta_{g(t)} + Z(t),$$

where $(Z(t))_{0 \leqslant t < T}$ is a $\mathrm{C}^{1,\infty}$-family of vector fields for $T \in (0, \infty]$.

We modify our standing assumption (A) in the following way.

**Assumption 4.1.** There are constants $\nu \geqslant 1$ and $\lambda \in \mathbb{R}$, such that it holds, for all $0 \leqslant s \leqslant T$,

$$\left(\mathsf{L}_s + \partial_s\right) r_s^2 \leqslant \nu + \lambda r_s^2, \tag{B}$$

Then theorems 3.2, 3.4 to 3.6 and 3.8 still hold true by replacing assumption (A) with assumption (B).

## 5 Applications and Examples

In [KP11, Theorem 1], Kuwada and Philipowski showed that $g(t)$-Brownian motion on $M$ cannot explode up to time $T$, if the family of Riemannian metrics evolves under a backwards super Ricci flow, i.e. $\partial_t g \leqslant \mathrm{Ric}$, generalising a fundamental principle that lower bounds on Ricci curvature imply the non-explosion of Brownian motion in the case of a static manifold, cf. [Tha16; HT94; Wan14].

In the evolving case, it is therefore natural to define

$$\mathscr{R}_Z^{g(t)}(X, Y) := \mathrm{Ric}^{g(t)}(X, Y) - \left\langle \nabla_X^{g(t)} Z(t), Y \right\rangle_{g(t)} - \frac{1}{2}\partial_t g(t)(X, Y),$$

for vector fields $X, Y \in \Gamma(\mathrm{T}M)$.

Following the argument in [Che17, p. 782] and the Proof of [KP11, Lemma 6] if there exist nonnegative functions $\varphi \in C[0, \infty)$ and $h \in C[0, T)$ such that $\mathscr{R}_Z^{g(t)} \geqslant -h(t)\varphi(r_t)$, we find a non-increasing function $F$ satisfying $\lim_{r \to 0} r\, F(r) < \infty$ such that

$$\left(\mathsf{L}_t + \partial_t\right) r_t(x) \leqslant F(r_t(x)) + h(t) \int_0^{r_t(x)} \varphi(s)\mathrm{d}s + |Z(t, o)|_t\,.$$

**Remark 5.1.** Assumption (B) can be ensured by the following two conditions:

(B1) There is a nonnegative function $k \in C[0, T)$ such that $\mathscr{R}_Z^{g(t)} \geqslant -k(t)$ for all $t \in [0, T)$.

(B2) There are two nonnegative functions $C_1, C_2 \in C[0, T)$ such that, for all $0 \leqslant t \leqslant T$, we have

$$\mathrm{Ric}_{g(t)} \geqslant -C_1(t)(1 + r_t^2) \qquad \text{and} \qquad \partial_t r_t + \left\langle Z(t), \nabla^{g(t)} r_t\right\rangle_{g(t)} \leqslant C_2(t)(1 + r_t).$$

By (B1), we can use the Laplace comparison theorem, so that

$$\begin{aligned}
\left(\mathsf{L}_s + \partial_s\right) r_s^2 &= 2r_s\left(\mathsf{L}_s + \partial_s\right) r_s + 2\left|\nabla^{g(t)} r_s\right|_{g(t)}^2 \\
&= 2r_s \Delta_{g(s)} r_s + 2(\partial_s + Z(s)) r_s + 2 \\
&\leqslant 2r_s \sqrt{(m-1)C_1(t)(1 + r_t^2)} \coth\left(\sqrt{(m-1)C_1(t)(1 + r_t^2)}\right) + 2C_2(t) r_t(1 + r_t^2) + 2,
\end{aligned}$$

where we used assumption (B2) in the third step. Using the inequality $\coth s \leqslant 1 + s^{-1}$, we see that (B) is satisfied.